\documentclass[twoside,12pt]{article}
\usepackage[all]{xy}
\usepackage{amssymb,amsmath,latexsym}

\usepackage{amsfonts}
\usepackage{ulem}
\usepackage{amscd}
\usepackage{amsxtra}
\usepackage{mathrsfs}

\newcommand{\field}[1]{\mathbb{#1}}

\newcommand{\R}{\field{R}}

\newcommand{\Z }{\field{Z}}

\newtheorem{theorem}{Theorem}[section]
\newtheorem{lemma}[theorem]{Lemma}
\newtheorem{proposition}[theorem]{Proposition}

\newtheorem{definition}[theorem]{Definition}
\newtheorem{corollary}[theorem]{Corollary}
\newtheorem{example}[theorem]{Example}
\newtheorem{remark}[theorem]{Remark}

\numberwithin{equation}{section}

\newcommand{\cqfd}
{\hspace{1cm}
\rule{2mm}{2mm}%
\medbreak%
\par%
}

\def\pr{{\parindent0pt {\bf Proof.\ }}}
\def\cqfd
{\hspace{1cm}
\rule{2mm}{2mm}%
\medbreak%
\par%
}

\def\tri{{\rm Tri}}

\def\Der{{\rm Der}}
\def\der{{\rm der}}

\def\GDer{{\rm  GDer}}

\def\Innder{{\rm Innder}}
\def\innder{{\rm innder}}
\def\InnBi{{\rm  InnBi}}

\def\Innbi{{\rm  Innbi}}
\def\InnGd{{\rm  InnGd}}
\def\End{{\rm  End}}

\def\Hom{{\rm  Hom}}

\def\Ker{{\rm Ker}}
\def\Im{{\rm Im}}

\author{}

\begin{document}
\title{Derivations and the first cohomology group of trivial extension algebras}

\date{}
\maketitle \vspace*{-1.5cm}
\thispagestyle{empty}
\begin{center}
\author{Driss Bennis$^{1}$  and Brahim Fahid$^{2}$}

\small{Department of Mathematics, Faculty of Sciences,  B.P. 1014, Mohammed V University in Rabat,   Morocco.}
\end{center}\vspace{-0,3cm}
\hspace{4cm} \small{ 1: d.bennis@fsr.ac.ma; driss$\_$bennis@hotmail.com}\\
$\mbox{}$\hspace{4cm} \small{ 2:  fahid.brahim@yahoo.fr}

\bigskip\bigskip
\noindent{\large\bf Abstract.}  In this paper we investigate in details derivations on trivial extension algebras. We obtain generalizations of both known results on derivations on triangular matrix algebras and a known result on first cohomology group of trivial extension algebras. As a consequence we get the characterization of trivial extension algebras on which every derivation is inner.  We show that, under some    conditions,  a trivial extension algebra on which every derivation is inner has necessarily a   triangular matrix representation. The paper starts with detailed study (with examples) of the relation between the trivial extension algebras and the  triangular matrix algebras.\bigskip

\small{\noindent{\bf 2010 Mathematics Subject Classification.}  16W25, 15A78,  	16E40.}\smallskip

\small{\noindent{\bf Key Words.}  trivial extension algebra;   triangular matrix algebra;  derivation; inner derivation; module generalized derivation;   central inner bimodule homomorphism; cohomology group.}\bigskip



\section{Introduction}
Throughout the paper $\mathcal{R}$ will denote a commutative ring with unity,  $\mathcal{A}$ will be a unital $\mathcal{R}$-algebra with center $Z(\mathcal{A})$ and $\mathcal{M}$ will be a unital $\mathcal{A}$-bimodule.\medskip

Recall that an $\mathcal{R}$-linear map  $\mathcal{D}$ from $\mathcal{A}$ into $\mathcal{M}$ is said to be a derivation if $  \mathcal{D}(ab)= \mathcal{D}(a)b+a\mathcal{D}(b)$ for all $a,b \in \mathcal{A}.$
It is known that the sum of two derivations in $ \mathcal{A}$  with values in  $ \mathcal{M}$ is also a derivation. This defines the structure of a group on the set of all derivations  in $ \mathcal{A}$  with values in  $ \mathcal{M}$ denoted by  $\Der(\mathcal{A},\mathcal{M})$. In particular, when  $\mathcal{M}=\mathcal{A}$, we simply set   $\Der(\mathcal{A}):=\Der(\mathcal{A},\mathcal{A})$. A derivation $\mathcal{D}\in \Der(\mathcal{A},\mathcal{M})$ is said to be inner if it is of the form  $\mathcal{D}(a)=[m_{0},a] $  for some $m_{0} \in \mathcal{M}$, where $[-,- ]$  stands for the Lie bracket.   Also, it is known that the set of all inner derivations  in $ \mathcal{A}$  with values in  $ \mathcal{M}$ is a subgroup of $\Der(\mathcal{A},\mathcal{M})$. It will be denoted by $\Innder(\mathcal{A},\mathcal{M})$, and when  $\mathcal{M}=\mathcal{A}$, we simply set   $\Innder(\mathcal{A}):=\Innder(\mathcal{A},\mathcal{A})$.  It is a well-known fact that a derivation needs not to be inner. Namely, the well-known first cohomology group    $ H^1(\mathcal{A}) := \Der(\mathcal{A}) /  \Innder(\mathcal{A})$ measures how much the group of all derivations on   $\mathcal{A}$ differs from the group  on inner derivations.\medskip

 Several authors have been interested in finding suitable conditions under which every derivation on a given algebra  is  inner.  In  \cite{S}  Coelho and  Milies proved that every derivation on the upper triangular matrices is inner. A similar claim holds for nest algebras in \cite{K}. In \cite{D} Benkovi\v{c}  proved that every derivation of the block upper triangular matrix algebra is   inner and every  derivation of    triangular matrix algebras, under some conditions, is inner  (see \cite{W} and \cite{FM}). One of the main results in this paper generalizes the above result to the context of trivial extension algebras (see Section 2 for details about this construction). It is worth noting, that the notion of trivial extension algebras, which is a  generalization of   triangular matrix  algebras, has proved to be an excellent construction for providing interesting examples (see for instance \cite{Z}). In recent years,  some results on mapping problems of   triangular matrix  algebras have been extended to trivial extensions (see for instance \cite{H, MP,A}). In this paper we follow this stream and investigate in details derivation on trivial extension algebras. Our study leads to generalizations of both known and recent results on the first cohomology group of both trivial extention algebras and   triangular matrix  algebras (see Section 4). As mentioned in the abstract, we then naturally get  a characterization of trivial extension algebras (and then   triangular matrix  algebras) on which every derivation is inner (see Corollaries \ref{cor-princ-inner-triv} and \ref{cor-princ-inner-matr}). However, we will show that, under some  conditions,  a trivial extension algebra on which every derivation is inner has   necessarily a   triangular matrix  representation.\medskip

The paper is organized in the following way: \medskip

In Section  2 we study  the relation between   trivial extension algebras and triangular matrix algebras. More precise, we give conditions under which a trivial extension algebra has a triangular matrix representation (see Proposition \ref{prop-triv-trian} and Corollaries  \ref{cor-2-triv-trian} and \ref{cor-3-triv-trian}). We also discuss some interesting properties of triangular matrix algebras that   trivial extension algebras could also  possess without necessarily having a   triangular matrix representation (see Example \ref{exm-principal} and Proposition \ref{lem-idem-equa}).\medskip

 In Section 3   we investigate derivations on trivial extensions. Description of the form of derivations on trivial extension algebras (see Lemma \ref{lem-der}) shows that, in order to well understand these derivations, the notion of a module generalized derivations should be studied (see the first paragraph of Section 3). The module generalized derivations are extensions of both bimodule homomorphisms and the classical ``ring" generalized derivations. The first part of Section 3 is devoted to the study of module generalized derivations. The notion of inner (resp.  central) module generalized derivations is introduced (see Definitions  \ref{def-inn-g-mod} and \ref{def-inn-bimo}). These two notions have played a principal role in our study and,  in fact, they can be considered as a new approach of the study of both derivations and first cohomology group of trivial extension algebras. So, a characterization of when  module generalized derivations and bimodule homomorphisms are inner are given (see Theorem \ref{thm-inn-bimo}). As a consequence we describe   the
form of inner derivations on trivial extension algebras  (Theorem \ref{thm-form2}). It is worth noting that our results are mainly inspired by the study of derivations on triangular matrix algebras (please, see Remark \ref{rem-deriv} and Example \ref{exm-inn-bimo} which discuss    the relation between  derivations    on trivial extension algebras and those on   triangular matrix  algebras).\medskip

Section 4  is devoted to the study of the first cohomology group of trivial extension algebras. In this context the notion of restricted first cohomology groups is introduced. Using   Proposition \ref{prop-triv-trian}, which relates the restricted first cohomology group with the classical first  cohomology group, one can show that the study of the first cohomology group of trivial extension algebras is based on the investigation of the restricted first cohomology groups. The main result in this section (Theorem \ref{thm-princ})   relates the restricted first cohomology groups with the first cohomology group of the base ring and the quotient group  of the group of all bimodule homomorphisms by the subgroup of all central inner bimodule homomorphisms. This leads to, Corollary \ref{cor-triv-princ}, a generalization of  both  the classical result \cite[Theorem 5.5]{CMRS} and the recent result  \cite[Theorem  4.4]{AGST} which uses a purely homological argument  (also, compare it with \cite[Theorem 2.5]{MP}). As a consequence we get   a characterization of trivial extension algebras  on which every derivation is inner (see Corollaries \ref{cor-princ-inner-triv} and \ref{cor-princ-inner-matr}). Several attempts were effectuated in order to give an example of a trivial extension algebra not having a triangular matrix
representation and on which every derivation is inner. This leads to a natural question: How do the innerness of derivations on trivial extension algebras affects its structure? (See Question posed at the end of the paper). In this context, Theorem \ref{thm-triv-is-trian} shows that, under some conditions, the innerness of all derivations on a trivial extension algebra  implies that it has a   triangular matrix representation.

\section{Trivial extension algebras and    triangular matrix algebras}

This section is devoted    to a discussion on the relation between trivial extension algebras  and     triangular matrix algebras. First we recall the definition of these classical constructions and we give some basic known results.\bigskip

As mentioned in the introduction, we use   $\mathcal{R}$ to denote a commutative ring with unity,  $\mathcal{A}$ is a unital $\mathcal{R}$-algebra with center $Z(\mathcal{A})$ and $\mathcal{M}$ is   a unital $\mathcal{A}$-bimodule.\medskip

Let $A$ and  $B$ be two  $\mathcal{R}$-algebras and let $M$ be an $(A,B)$-bimodule.
The  set   $$ \tri(A;M;B)=\{\left(                                                                                                                                      \begin{array}{cc}                                                                                                                                      a & m \\                                                                                                                                0 & b                                                                                                                               \end{array}                                                                                                                              \right) \mid a\in A, m\in M, b\in B\}$$
endowed  with the usual matrix operations is an $\mathcal{R}$-algebra called a   triangular matrix $\mathcal{R}$-algebra  (see \cite{GJS} and \cite{W} for more details about this construction). Following  \cite[Chapter 5]{GJS}, an algebra  $S$ is said to have a   triangular matrix representation if $S$ is   isomorphic to a triangular matrix algebra. By \cite[Theorem 5.1.4]{GJS},  a unital algebra $\mathcal{A}$ has a   triangular matrix representation  if  there exists a non trivial idempotent $ e \in \mathcal{A}$ such that $ (1-e)\mathcal{A}e=0 $. Namely, in this case, $\mathcal{A}$ is isomorphic to $\tri(e\mathcal{A}e; e\mathcal{A}(1-e); (1-e)\mathcal{A}(1-e))$. Clearly, a   triangular matrix  algebra is an example of a noncommutative algebra; its  center $ Z(\tri(A;M;B))$ is  the set of elements  $\left(                                                                                                                                      \begin{array}{cc}                                                                                                                                      a & 0 \\                                                                                                                                0 & b                                                                                                                               \end{array}                                                                                                                              \right)$, where  $a\in Z(A)$ and   $b\in  Z(B) $, such that   $am=mb$ for  all $m\in M$.  As interesting examples of   triangular matrix algebras one can cite the  (classical) upper triangular matrix algebras, the block upper triangular matrix algebras, one-point extension algebras and the nest algebras (see,  for instance, \cite{D}). \bigskip

As shown below, a triangular algebra can be obtained as a special case of  trivial extension algebra. Recall that the direct product $\mathcal{A}\times \mathcal{M}$ together with the addition pairwise  scalar product and the algebra multiplication defined by  $ (a,m)(b,n)=(ab,an+mb)$ for all $a,b\in \mathcal{A}$ and $m,n \in \mathcal{M}$, is a unital algebra which is called a trivial extension of $\mathcal{A}$ by $\mathcal{M}$ and will be denoted by $\mathcal{A}\ltimes \mathcal{M}$. \\
The center of a trivial extension algebra $\mathcal{A}\ltimes \mathcal{M}$ is determined as follows:
\begin{equation*}
Z(\mathcal{A}\ltimes \mathcal{M})=\{(a,m)\mid a\in Z(\mathcal{A})\; \mathrm{and}\; [b,m]=0=[a,y]\;  \mathrm{for} \; \mathrm{all} \;  b\in \mathcal{A} \; \mathrm{and}\;  y\in  \mathcal{M}\}.
\end{equation*}
Note also that $Z(\mathcal{A}\ltimes \mathcal{M})= \pi_{\mathcal{A}}(Z(\mathcal{A}\ltimes \mathcal{M}))\times \pi_{\mathcal{M}}(Z(\mathcal{A}\ltimes \mathcal{M}))$,
where $ \pi_{\mathcal{A}} : \mathcal{A}\ltimes \mathcal{M} \longrightarrow \mathcal{A}$ and $\pi_{\mathcal{M}} : \mathcal{A}\ltimes \mathcal{M} \longrightarrow \mathcal{M}$ are the natural projections given by $\pi_{\mathcal{A}}(a,m)=a$ and $\pi_{\mathcal{A}}(a,m)=m$  for all $(a,m)\in \mathcal{A}\ltimes \mathcal{M}  $ (see \cite{A}). \bigskip

It is well-known that every triangular matrix algebra can be viewed as   a trivial extension algebra. Indeed,  $ \tri(A;M;B)$ is isomorphic to $( A  \times B)\ltimes M$ as an $\mathcal{R}$-algebra where    $M$ is viewed as an $A\times B$-bimodule via the module actions given by $(a,b)m=am$ and $m(a,b)=mb$ for all $(a,b) \in A\times B$ and  $ m\in M$. However, a  trivial extension algebra has not necessarily a   triangular matrix representation. To give an appropriate example, we set the following obvious but important result.

\begin{proposition}\label{prop-triv-trian}
A trivial extension algebra  $\mathcal{A}\ltimes \mathcal{M}$  has  a   triangular matrix representation if and only if there exists a non trivial idempotent $e$ of $\mathcal{A}$ such that $(1-e)\mathcal{A}e=0$ and $(1-e)\mathcal{M}e=0$. In particular, if  $\mathcal{A}\ltimes \mathcal{M}$  has  a triangular matrix representation, then so has  $\mathcal{A}$.
\end{proposition}
\pr If  $\mathcal{A}\ltimes \mathcal{M}$  has  a   triangular matrix representation, then  there exists a non trivial idempotent $ E=(e,m)\in \mathcal{A}\ltimes \mathcal{M}$ such that $ (1-E) \mathcal{A}\ltimes \mathcal{M}E=0 $. Clearly, $e$ is  a non trivial idempotent of
$\mathcal{A}$ which satisfies  $(1-e)\mathcal{A}e=0$ and $(1-e)\mathcal{M}e=0$.\\
The converse holds by considering the non trivial idempotent $(e,0)$ of $\mathcal{A}\ltimes \mathcal{M}$.\cqfd

As noted above the decision of whether an algebra $S$  has  a   triangular matrix representation depends of the existence of an appropriate idempotent. Moreover, if such an  idempotent exists,  say  $e$ ($fSe=0$, where $f=1-e$), then, for  $M:=eSf$,   $emf=m$ for every $m\in M$. This property on $\mathcal{M}$ and its consequences (see Lemma \ref{lem-idem-equa} below) play a crucial role in proving some interesting results (see, for instance, \cite{D, H, A}). In \cite[Example 3.13]{H} and \cite[Example 2.6]{A}, examples of  trivial extension algebras $\mathcal{A}\ltimes \mathcal{M}$ with the suitable idempotent of $\mathcal{A}$ without having  a   triangular matrix representation are given.  A deep observation of these  examples leads to consider  the following particular case of Proposition \ref{prop-triv-trian}.\\

In the following result we use the following well-known fact: If $\mathcal{N}$ is an  $\mathcal{A}$-bimodule, then it can be used to define an $\mathcal{A}\ltimes \mathcal{M}$-bimodule
 via the module actions given by $(a,m)n=an$ and $n(a,m)=na$ for all $(a,m) \in \mathcal{A}\ltimes \mathcal{M}$ and $n\in \mathcal{N}$.

\begin{corollary}\label{cor-2-triv-trian} For an   $\mathcal{A}$-bimodule $\mathcal{N}$, the following assertions are equivalent.
\begin{enumerate}
    \item The trivial extension algebra  $(\mathcal{A}\ltimes \mathcal{M})\ltimes\mathcal{N} $  has a   triangular matrix representation.
    \item  There exists a non trivial idempotent $e$ of $\mathcal{A}$ such that $(1-e)\mathcal{A}e=0$, $(1-e)\mathcal{M}e=0$ and $(1-e)\mathcal{N}e=0$.
\end{enumerate}
\end{corollary}

In particular, when $\mathcal{A}\ltimes \mathcal{M}$ is a   triangular matrix  algebra $\tri(A;M;B)$, we get the following result.

\begin{corollary}\label{cor-3-triv-trian} Let $A$ and  $B$ be two  $\mathcal{R}$-algebras and $M$ be an $(A,B)$-bimodule.  Let   $\mathcal{S}= \tri(A;M;B)$ and  $\mathcal{N}$ be an $\mathcal{S}$-bimodule. The following assertions are equivalent.
\begin{enumerate}
    \item The trivial extension algebra  $\mathcal{S}\ltimes\mathcal{N} $  has a   triangular matrix representation.
    \item  There exists a non trivial idempotent $e=(e_a,e_b)\in A\times B$ such that $(1-e_a)Ae_a=0$, $(1-e_b)Be_b=0$, $(1-e_a)M e_b=0$ and $(1-e)\mathcal{N}e=0$.
\end{enumerate}
In particular, if in addition $A$ and $B$ have only trivial idempotents, then the idempotent $e$ (in the second assertion) is necessarily $(1,0)$.
\end{corollary}

Now we can  easily construct the desired example (compare with  \cite[Example 3.13]{H} and \cite[Example 2.6]{A}).

\begin{example}\label{exm-principal}  Let $A$ and  $B$ be two  $\mathcal{R}$-algebras which have   only trivial idempotents. Let $M$ be both  an $(A,B)$-bimodule  and a  $(B,A)$-bimodule.  Let   $\mathcal{S}:= \tri(A;M;B)$.  We make $M$ into an  $ \mathcal{S}$-bimodule by defining    $((a,b),m)m':= bm'$  and $m'((a,b),m):= m'a$ for all $(a,b) \in A\times B$ and $m,m'\in  M$.  As an $ \mathcal{S}$-bimodule, $M$ will be denoted as $ \mathcal{N}$.\\
   Then, using Corollary \ref{cor-3-triv-trian}, one can show that the trivial extension of $\mathcal{S}$ by $\mathcal{N}$ does not have a   triangular matrix representation, because the idempotent $F=((0,1),0)$ of $\mathcal{S}$ satisfies   $Fn  (1-F)=n$ for every $n\in \mathcal{N}$.
\end{example}

The last condition is equivalent to some other interesting conditions  as shown in the following result which  can be proven easily   (compare with the remark given before \cite[Theorem 2.2]{A}).

\begin{proposition}\label{lem-idem-equa}  Consider a non-trivial idempotent $e$ of an algebra $S$ and set $f=1-e$. For every $S$-bimodule $N$, the following assertions are equivalent.
\begin{enumerate}
    \item For every $m\in N$, $emf=m$.
    \item For every $m\in N$, $fm =0=me$.
    \item For every $m\in N$, $em =m=mf$.
  \item For every $m\in N$ and $a\in S$, $am=eaem$ and $ma=mfaf$.
\end{enumerate}
\end{proposition}

  In this paper, we sometimes consider a trivial extension $\mathcal{A}\ltimes \mathcal{M}$ which satisfies the following property (satisfied by the ring of Example \ref{exm-principal}):  There is a non trivial idempotent $e$ of  $\mathcal{A}$  such that $e\mathcal{A}f=0$ (where  $f=1-e$) and  $emf=m$ for every $m\in \mathcal{M}$ (thus  $\mathcal{A}\ltimes \mathcal{M}$  satisfies all the equivalent properties of Proposition \ref{lem-idem-equa}). We refer to such an algebra by a trivial extension of type $(\star)$ (this trivial extension algebras were first considered in \cite{A}). Note that in this case, we have $(e,0)(a,m)(e,0)=(eae,0)$ and $(f,0)(a,m)(f,0)=(faf,0)$. Then, $(e,0)(\mathcal{A}\ltimes \mathcal{M})(e,0)\cong  e\mathcal{A} e $ and $(f,0)(\mathcal{A}\ltimes \mathcal{M})(f,0)\cong f \mathcal{A} f$. And, using assertion (4) of Proposition \ref{lem-idem-equa},  one can see that $\mathcal{M}$ is  a left $e\mathcal{A} e $-module and a right $f \mathcal{A} f$-module via the module actions  $am=eaem$ and $ma=mfaf$ for all $m\in \mathcal{M}$ and $a\in \mathcal{A}$.\medskip


\section{Derivations on trivial extension algebras}
In this section we investigate derivations on trivial extension algebras.  \medskip

We begin with \cite[Lemma 2.1]{A} which describes derivations on trivial extensions (see also \cite[Proposition 2.2]{MP}).   We set \cite[Lemma 2.1]{A} using the following  terminology:\medskip

 The set  of all  $\mathcal{A}$-bimodule homomorphisms $ f: \mathcal{M}\longrightarrow \mathcal{N} $ is an  $\mathcal{R}$-module  denoted by  $\Hom_{\mathcal{A}-\mathcal{A}}(\mathcal{M},\mathcal{N})$ and,  when $ \mathcal{M}=\mathcal{N}$, it   is denoted by $ \End_{\mathcal{A}-\mathcal{A}}(\mathcal{M})$.

 Following the notation in \cite{SM},   let $\mathcal{E}( \mathcal{M})$ denote  the $\mathcal{R}$-submodule of $ \Hom_{\mathcal{A}-\mathcal{A}}( \mathcal{M} ,\mathcal{A})$ consisting of all $\mathcal{A}$-bimodule homomorphisms  $f$  such that $f(m)n+mf(n) = 0$ for all $m,\,n\in  \mathcal{M}$.  A linear map $ S: \mathcal{M}\longrightarrow \mathcal{M} $ is called a \textit{module generalized $d$-derivation}  for some $ d\in \Der(\mathcal{A})$ if  $ S(am)=aS(m)+d(a)m$ and $S(ma)=S(m)a+md(a)$  for  all $m\in \mathcal{M} $ and $a\in \mathcal{A}$. If   there is no ambiguity about the associated derivation $d$, $S$ will be simply called a module generalized derivation. Clearly the set of all module generalized derivations $ S: \mathcal{M}\longrightarrow \mathcal{M} $ is a group denoted   by $  \GDer(\mathcal{M})$. Compare this notion with the notion of the generalized derivation on  modules introduced in \cite{GMA} and the classical generalized derivation on rings \cite{Bre91}. Also, it is clear that $\End_{\mathcal{A}-\mathcal{A}}(\mathcal{M})$ is a subgroup of $  \GDer(\mathcal{M})$. Namely, every $\mathcal{A}$-bimodule homomorphism can be considered as a module generalized derivation associated to the zero derivation.

\begin{lemma}[\cite{A}, Lemma 2.1]\label{lem-der}
 Every linear map $ \mathcal{D}: \mathcal{A}\ltimes \mathcal{M} \longrightarrow \mathcal{A}\ltimes \mathcal{M}$ has the form
  \begin{equation*}
\mathcal{D}(a,m)=(\mathcal{D}_{\mathcal{A}}(a)+T(m),\mathcal{D}_{\mathcal{M}}(a)+S(m))\;\;\; (a\in \mathcal{A},\; m\in \mathcal{M})
  \end{equation*}
  for some linear maps $\mathcal{D}_{\mathcal{A}}: \mathcal{A}\longrightarrow \mathcal{A}, \mathcal{D}_{\mathcal{M}}: \mathcal{A}\longrightarrow \mathcal{M}, T: \mathcal{M}\longrightarrow \mathcal{A}$ and $S:\mathcal{ M}\longrightarrow \mathcal{M}$. We will simply write $\mathcal{D}=(\mathcal{D}_{\mathcal{A}}+T,\mathcal{D}_{\mathcal{M}}+S )$.\\
 Moreover, $\mathcal{D}$ is a derivation if and only if  $\mathcal{D}_{\mathcal{A}}$ and $\mathcal{D}_{\mathcal{M}}$ are derivations, $  T\in \mathcal{E}( \mathcal{M})$  and $S$ is a module generalized $\mathcal{D}_{\mathcal{A}}$-derivation.
\end{lemma}

Our aim in this section is both, the study of the relation  between   module generalized derivations and their associated derivations, and the innerness of derivations on trivial extension algebras.\medskip

Let us start with the following remark which sheds light on some significant differences  between the forms of derivations on trivial extension algebras and of those on triangular matrix algebras.

\begin{remark}\label{rem-deriv}
 \textnormal{
  \begin{itemize}
     \item [(i)]   Note that when $ \mathcal{A}\ltimes \mathcal{M}=\tri(A,\mathcal{M},B)$ (where $\mathcal{A}=A\times B $), the derivation $\mathcal{D}_{\mathcal{M}}$ is necessarily inner (see \cite[Theorem 2.2.1]{W} and \cite{FM}).  However, for a trivial extension algebra which has not a  triangular matrix representation  this does not hold true in general. To see this, it is sufficient to consider a derivation $ d:\mathcal{A}\longrightarrow \mathcal{M}$ which is not inner and consider $\mathcal{D}=(0,d):\mathcal{A}\ltimes \mathcal{M}\longrightarrow \mathcal{A}\ltimes \mathcal{M}$ which is a derivation on $\mathcal{A}\ltimes \mathcal{M}$.\\
Nevertheless, for a trivial extension of type $(\star)$,  the derivation $\mathcal{D}_{\mathcal{M}}$ is necessarily inner. Indeed, for $a\in \mathcal{A}$, we have 
\begin{eqnarray*}
     \mathcal{D}_{\mathcal{M}}(a) &=&\mathcal{D}_{\mathcal{M}}(eae+fae+faf) \\
            &=&  \mathcal{D}_{\mathcal{M}}(eae) +\mathcal{D}_{\mathcal{M}}(fae)+\mathcal{D}_{\mathcal{M}}(faf)\\
     &=&  ea \mathcal{D}_{\mathcal{M}}(e) + \mathcal{D}_{\mathcal{M}}(f)af\\
  &=&   a \mathcal{D}_{\mathcal{M}}(e) + \mathcal{D}_{\mathcal{M}}(f)a \\
  &=&   a m - m a,
         \end{eqnarray*}
where $m=\mathcal{D}_{\mathcal{M}}(e)=-\mathcal{D}_{\mathcal{M}}(f)$.
     \item[(ii)]  Also,  when $ \mathcal{A}\ltimes \mathcal{M}=\tri(A,\mathcal{M},B)$ (where $\mathcal{A}=A\times B $), the map  $T$ is zero (see \cite[Theorem 2.2.1]{W} and \cite{FM}).
In fact, this holds for every trivial extension of type $(\star)$. Indeed, for $m\in \mathcal{M}$, we have $$T(m)=T(em)=e\, T(m) \quad  \mathrm{and} \quad T(m)=T( mf)= T(m)\, f.$$  Therefore, $$T(m)=e\, T(m)=   e\,  T(m)\, f=0 .$$
However, the following example shows that the map $T$ needs  not to be zero in general.\\
Consider the trivial extension $ M_2(\Z/2\Z)\ltimes M_2(\Z/2\Z)$ where $ M_2(\Z/2\Z) $ is the algebra of $2\times 2$ matrices with entries from $\Z/2\Z$. Now consider the identity map $T:  M_2(\Z/2\Z)\to  M_2(\Z/2\Z)$. Since $M_2(\Z/2\Z)$ has characteristic $2$,  $T\in   \mathcal{E}( M_2(\Z/2\Z))$. Then, the linear map $ \mathcal{D}:   M_2(\Z/2\Z)\ltimes M_2(\Z/2\Z)   \to     M_2(\Z/2\Z)\ltimes M_2(\Z/2\Z)$ defined by $\mathcal{D}((a,b))=(T(b),0) $  for all $ a,b\in   M_2(\Z/2\Z)$ is   a  derivation with  $ T \neq 0 $.
   \end{itemize}
}
\end{remark}

In  \cite[Corollary 2.2.2]{W},  faithfulness of $\mathcal{M}$ was used to  show that a derivation $    \mathcal{D} =(\mathcal{D}_{\mathcal{A}}   ,\mathcal{D}_{\mathcal{M}} +S )$ on a   triangular matrix algebra $ \mathcal{A}\ltimes \mathcal{M}=\tri(A,\mathcal{M},B)$ (where $\mathcal{A}=A\times B $) is uniquely determined by $S$ and $\mathcal{D}_{\mathcal{M}}$.   Here we use a condition in terms of annihilators. Recall that the left annihilator, denoted $l.Ann_{\mathcal{A}}(\mathcal{M})$, of $\mathcal{M}$ is the set of all elements $r$ in $\mathcal{A}$ such that $r\mathcal{M}=0$. Similarly   the right annihilator,   $r.Ann_{\mathcal{A}}(\mathcal{M})$,   of $\mathcal{M}$ is defined. When $ \mathcal{A}\ltimes \mathcal{M}=\tri(A,\mathcal{M},B)$ (where $\mathcal{A}=A\times B $),  one can show easily that the condition  $l.Ann_{\mathcal{A}}(\mathcal{M})\cap r.Ann_{\mathcal{A}}(\mathcal{M})=0$ is equivalent to the fact that $\mathcal{M}$ is faithful as both a left $A$-module and a right  $B$-module. However, it is worth noting that for trivial extension algebras $l.Ann_{\mathcal{A}}(\mathcal{M})\cap r.Ann_{\mathcal{A}}(\mathcal{M})=0$ does not imply that $\mathcal{M}$ is a faithful $\mathcal{A}$-bimodule.\medskip

From the proof of \cite[Corollary 2.2.2]{W}, we can deduce the following result:\smallskip

\indent   Let $A$ and $B$ be unital algebras over a commutative ring $ \mathcal{R} $, and let $\mathcal{M}$ be a unital $(A,B)$-bimodule, which is faithful as both a left $A$-module and  a right $B$-module. Consider  a linear map $ S: \mathcal{M}\longrightarrow \mathcal{M} $. If  there are linear maps $ d_{A} : A\longrightarrow A $ and $ d_{B} : B\longrightarrow B $ such that  $S(am)=aS(m)+d_{A}(a)m$ and $S(ma)=S(m)a+md_B(a)$  for  all $m\in \mathcal{M} $, $a\in A$ and $b\in B$, then $ d_{A}$ and $ d_{B}$ are derivations.\\

 Using similar arguments as used in the proof of \cite[Corollary 2.2.2]{W},  we get the following extension to the context of  trivial extension algebras.

\begin{proposition} \label{prop-der-faith} Assume that $l.Ann_{\mathcal{A}}(\mathcal{M})\cap r.Ann_{\mathcal{A}}(\mathcal{M})=0$. Consider  a linear map $ S: \mathcal{M}\longrightarrow \mathcal{M} $. If  there is a linear map $ d  : \mathcal{A}\longrightarrow  \mathcal{A} $   such that  $S(am)=aS(m)+d (a)m$ and $S(ma)=S(m)a+md (a)$  for  all $m\in \mathcal{M} $ and $a\in A$, then $ d$ is a derivation. Consequently, $S$ is  a  module generalized $d$-derivation.
\end{proposition}

Now we investigate  a particular case of  module generalized derivations. Let us first give the following observation.\medskip

Recall that when $ \mathcal{A}\ltimes \mathcal{M}=\tri(A,\mathcal{M},B)$ (where $\mathcal{A}=A\times B)$, the homomorphism  $ S:\mathcal{ M}\longrightarrow \mathcal{M}$ associated to an inner derivation is of the form $ S(m)=a_{0}m-mb_{0}$ for fixed $ (a_{0}, b_{0})\in  A\times B $ and  for all $ m\in \mathcal{M}$ (see \cite[Proposition 2.2.3]{W}, \cite[Proposition 3.3]{D}   and  \cite{GMA}).   Note that   such a homomorphism   can be written in the form similar to that of the classical inner derivation. Indeed, for every $ m\in \mathcal{M}$, $$S(m)=a_{0}m-m b_{0}=(a_0,b_0) m-m(a_0,b_0)=  c_0m-mc_0,$$
where $c_0 =(a_0,b_0)$.\\
In general, for every $c\in   \mathcal{A}$, the linear map  $ S:\mathcal{ M}\longrightarrow \mathcal{M}$, defined by $S(m) =  c m-mc$ for all $ m\in \mathcal{M}$, is a module generalized derivation associated to the inner derivation $d:\mathcal{ A}\longrightarrow \mathcal{A}$ defined by $d(a) =ca-ac$ for all $ a\in \mathcal{A}$. This observation leads us to introduce the following notion which plays an important role in the sequel. 

\begin{definition}\label{def-inn-g-mod}
 \textnormal{A  module generalized derivation   $ S:\mathcal{ M}\longrightarrow \mathcal{M}$  is said to be inner if there exists  $ a_{0}\in \mathcal{A}$ such that $ S(m)=a_{0}m-ma_{0}$ for all $ m\in \mathcal{M}$.
}
\end{definition}

\noindent\textbf{Notation.} 
The expression $ a_{0}m-ma_{0}$ in Definition \ref{def-inn-g-mod} above will be also noted by $[a_{0},m]$. Namely, we will use the Lie bracket in the following three contexts:
\begin{itemize}
  \item For an inner derivation  on $\mathcal{A}$ associated to an element $a_{0}\in \mathcal{A}$:
$$\begin{array}{rcll}
     [a_{0},-]:  &  \mathcal{A}&\longrightarrow   &  \mathcal{A} \\
   &  a &\longmapsto &  [a_{0},a]=a_{0}a-aa_{0} 
\end{array}
$$
\item For an inner derivation  from $\mathcal{A}$ to $\mathcal{M}$  associated to an element $m_{0}\in \mathcal{M}$: 
$$\begin{array}{rcll}
     [m_{0},-]:  &  \mathcal{A}&\longrightarrow   &  \mathcal{M} \\
   &  a &\longmapsto &  [m_{0},a]=m_{0}a-am_{0} 
\end{array}
$$
  \item For an inner (module generalized) derivation  from $\mathcal{M}$ to $\mathcal{M}$  associated to an element $a_{0}\in \mathcal{A}$: 
$$\begin{array}{rcll}
     [a_{0},-]:  &  \mathcal{M}&\longrightarrow   &  \mathcal{M} \\
   &  m &\longmapsto &  [a_{0},m]=a_{0}m-ma_{0} 
\end{array}
$$
 \end{itemize}

\begin{remark}\label{rem-def-innBi}
 It is evident that the set of all inner module generalized derivations $\mathcal{M}\longrightarrow \mathcal{M}$  is a group.   It will be denoted by $\InnGd_{\mathcal{A}}(\mathcal{M})$.\\
 Let us also denote  $\InnBi_{\mathcal{A}}(\mathcal{M}):=  \InnGd_{\mathcal{A}}(\mathcal{M})\cap \End_{\mathcal{A}-\mathcal{A}}(\mathcal{M}).$  The elements  of $\InnBi_{\mathcal{A}}(\mathcal{M})$ will be called inner $\mathcal{A}$-bimodule homomorphisms.  One can show that the group $\InnBi_{\mathcal{A}}(\mathcal{M})$ coincide with the set of all inner  module generalized derivations $[a_0,-]$, where $a_0\in\mathcal{A}$ satisfies $[a_0,a]\in l.Ann_{\mathcal{A}}(\mathcal{M})\cap r.Ann_{\mathcal{A}}(\mathcal{M})$ for every $a\in   \mathcal{A}$.
\end{remark}

It is also clear that, if    $ c\in Z(\mathcal{A})$,  the linear map  $ S:\mathcal{ M}\longrightarrow \mathcal{M}$, defined by $S(m) =  c m-mc$ for all $ m\in \mathcal{M}$, is an inner bimodule homomorphism. This kind of inner bimodule homomorphisms will be of particular interest. Namely, in the study of the first cohomology group of trivial extension algebras (see Section 4). For this reason, we introduce  the following notion. We use the terminology of \cite{FM} (see \cite[Remark 2.3(ii) and Definition 2.5]{FM}).  

\begin{definition}\label{def-inn-bimo}
 \textnormal{
  A  bimodule homomorphism   $ S:\mathcal{ M}\longrightarrow \mathcal{M}$  is said to be a central  inner  bimodule homomorphism   if there exists  $ a_{0}\in Z( \mathcal{A})$ such that $ S(m)=a_{0}m-ma_{0}$ for all $ m\in \mathcal{M}$.
}
\end{definition}

When $ \mathcal{A}\ltimes \mathcal{M}=\tri(A,\mathcal{M},B)$ (where $\mathcal{A}=A\times B)$, the central  inner  bimodule homomorphisms are called bimodule homomorphisms of the standard form (see \cite{D}).  These are exactly  the homomorphisms  $ S:\mathcal{ M}\longrightarrow \mathcal{M}$ defined by $S(m)=a_{0}m-mb_{0}$ for fixed $ (a_{0}, b_{0})\in  Z(A)\times Z( B)$ and for all $ m\in \mathcal{M}$.

\begin{remark}\label{rem-innBi2}
\begin{enumerate}
    \item  It is evident that the set of of all central inner  bimodule homomorphisms  $\mathcal{M}\longrightarrow \mathcal{M}$ is a subgroup of  $\InnBi_{\mathcal{A}}(\mathcal{M})$. It will be denoted by $\Innbi_{\mathcal{A}}(\mathcal{M})$.
\item  It is worth noting that    the inclusion $\Innbi_{\mathcal{A}}(\mathcal{M}) \subseteq \InnBi_{\mathcal{A}}(\mathcal{M})$ can be strict (see assertion (2) in Example \ref{exm-inn-bimo}). However, in the next result (Theorem \ref{thm-inn-bimo}), we show  that, when $l.Ann_{\mathcal{A}}(\mathcal{M})\cap r.Ann_{\mathcal{A}}(\mathcal{M})=0$,  the two groups  $\Innbi_{\mathcal{A}}(\mathcal{M})$ and $ \InnBi_{\mathcal{A}}(\mathcal{M})$  coincide.
\item  Note also that $\Innbi_{\mathcal{A}}(\mathcal{M})=0$ if and only if the center $Z(\mathcal{A})$   of $\mathcal{A}$ has a symmetric action on $\mathcal{M}$  (that is, $am=ma$ for every  $a\in Z(\mathcal{A})$ and $m\in \mathcal{M}$).
\end{enumerate}
\end{remark}
 
The next   result  gathers relations that exist between inner   module generalized  derivations and the associated (ring) derivations under the condition ``$l.Ann_{\mathcal{A}}(\mathcal{M})\cap r.Ann_{\mathcal{A}}(\mathcal{M})=0$". More relations will be given in Proposition \ref{prop-gen-bimo} without this condition.

\begin{theorem}  \label{thm-inn-bimo}   Assume that $l.Ann_{\mathcal{A}}(\mathcal{M})\cap r.Ann_{\mathcal{A}}(\mathcal{M})=0$. Consider a module generalized $d$-derivation $ S: \mathcal{M}\longrightarrow \mathcal{M} $, where $ d  : \mathcal{A}\longrightarrow \mathcal{A} $  is a derivation. 
\begin{enumerate}
    \item  If $S$ is inner, then $d$ is inner. Namely, if  for some $ a_{0}\in \mathcal{A}$,  $ S(m)= [a_0,m]$ for all $ m\in \mathcal{M}$, then $d(a)= [a_0,a]$ for all $ a\in \mathcal{A}$.
 \item    If $S$ is a bimodule homomorphism, then $d=0$. If in addition $S=[a_0,-]$ is inner, for some $ a_{0}\in \mathcal{A}$,    then  $a_0$ is in the  center $Z(\mathcal{A})$ of $ \mathcal{A}$ (i.e., $S$ is  a central inner  bimodule homomorphism).
\end{enumerate}
\end{theorem}  
\pr 1. Let $ a_{0}\in \mathcal{A}$ be such that $S(m)=[a_{0},m]$ for all $ m\in \mathcal{M}$. Using the fact that  $ S(am)=aS(m)+d(a)m$ and $S(ma)=S(m)a+md(a)$  for  all $m\in \mathcal{M} $ and $a\in \mathcal{A}$, we get, for all  $m\in \mathcal{M} $ and $a\in \mathcal{A}$, 
$$([a_{0},a]-d(a))m=0\quad \mathrm{and}\quad m([a_{0},a]-d(a))=0.$$
Then, by hypothesis, $d(a)=[a_{0},a]$, as desired.\smallskip

2.  If $S$ is  a bimodule homomorphism, then $d(a)m=0=md(a)$ for all $a\in \mathcal{A}$ and $m\in \mathcal{M}$. This shows that $d=0$.\\
Now let $ a_{0}\in \mathcal{A}$ be such that $S(m)=[a_{0},m]$ for all $ m\in \mathcal{M}$. Using the fact that  $ S(am)=aS(m) $ and $S(ma)=S(m)a$  for  all $m\in \mathcal{M} $ and $a\in \mathcal{A}$, we get, for all  $m\in \mathcal{M} $ and $a\in \mathcal{A}$
$$(a_{0}a -aa_0)m=0\quad \mathrm{and}\quad m(a_{0}a -aa_0)=0.$$
Then, by hypothesis, $ a_{0}a =aa_0 $, as desired.\cqfd

Next we give an example of  a non trivial inner  bimodule homomorphism and another example showing that, in order to get assertion (2) of Theorem \ref{thm-inn-bimo} above, the condition $l.Ann_{\mathcal{A}}(\mathcal{M})\cap r.Ann_{\mathcal{A}}(\mathcal{M})=0$ cannot be dropped. It is based on the following observation.

\begin{lemma} \label{lem-InnDer}   Every module generalized derivation  $ S: \mathcal{M}\longrightarrow \mathcal{M} $ is a bimodule homomorphism if the image of every derivation on   $ \mathcal{A}$ is in $l.Ann_{\mathcal{A}}(\mathcal{M})\cap r.Ann_{\mathcal{A}}(\mathcal{M})$.\\
  In other words, $  \GDer(\mathcal{M})=\End_{\mathcal{A}-\mathcal{A}}(\mathcal{M})$ if, for every $d\in \Der(\mathcal{A})$, $\Im(d)\subseteq l.Ann_{\mathcal{A}}(\mathcal{M})\cap r.Ann_{\mathcal{A}}(\mathcal{M})$, where   $\Im(d)$ denotes the image of $d$. 
\end{lemma} 
\pr Straightforward.\cqfd

\begin{example}\label{exm-inn-bimo}
Consider $\mathcal{R}$, $A$, $B$, $M$,  $\mathcal{S}$  and $\mathcal{N}$   as  in example \ref{exm-principal}. Assume that, as sets, $\mathcal{R}=A=B=M$. Then, the following assertions hold.
\begin{enumerate}
    \item Consider $M$ as an $(A,B)$-bimodule. Then, every  bimodule homomorphism  $ S:  M \longrightarrow M $ is inner. Note that $l.Ann_{A\times B}(M)\cap r.Ann_{A\times B}(M)=0$.
    \item  Every generalized derivation $ S:   \mathcal{N}  \longrightarrow  \mathcal{N}$ is an inner  generalized $d$-derivation for every derivation $d$ on  $\mathcal{S}$. Moreover, it is also a bimodule homomorphism.\\
Note that  $l.Ann_{\mathcal{S}}(\mathcal{N})\cap r.Ann_{\mathcal{S}}(\mathcal{N})=                                                                                                                                        \begin{pmatrix}                                                                                                                                    0 &  M\\                                                                                                                                0 & 0                                                                                                                              \end{pmatrix}  $.
\end{enumerate}
\end{example}
\pr 1. Consider  a bimodule homomorphism  $ S:  M \longrightarrow M $ and let $x\in M$. Then, $$S(x)=S((x,0)1)= (x,0) S(1)=xS(1)=(S(1),0)x=(S(1),0)x-x(S(1),0).$$
This shows that $S$ is inner.\\
2. First it is clear that every derivation on $\mathcal{S}$ is inner. Thus,  consider a derivation $d=[a ,-]$ on $\mathcal{S}$ for some $a =                                                                                                                                    \begin{pmatrix}                                                                                                                                      x_a & m_a\\                                                                                                                                0 & y_a                                                                                                                               \end{pmatrix}  $ in $\mathcal{S}$. Then, for every element
 $b =                                                                                                                                       \begin{pmatrix}                                                                                                                                       x_b & m_b\\                                                                                                                                0 & y_b                                                                                                                               \end{pmatrix}                                                                                                                               $ in $\mathcal{S}$,
 $$d(b)=[a ,b]=                                                                                                                                    \begin{pmatrix}                                                                                                                                        0&x_a m_b  +m_ay_b-x_bm_a-m_by_a \\                                                                                                                                0 & 0                                                                                                                               \end{pmatrix}                                                                                                                                $$
This shows that $d(b)n=0=nd(b)$ for all $n\in \mathcal{N}$. Then, every  generalized derivation $ S:   \mathcal{N}  \longrightarrow  \mathcal{N}$ is a bimodule homomorphism. It remains to prove that every  bimodule homomorphism is inner. Then consider a bimodule homomorphism  $ S:  \mathcal{N} \longrightarrow \mathcal{N}  $ and let $x\in \mathcal{N}$. Then,  $S(x)=S(((0,x),0)1)= ((0,x),0)S(1)=xS(1)=((0,S(1)),0) x=((0,S(1)),0)x-x((0,S(1)),0).$
Therefore, $S$ is inner.\cqfd

The following result, which is a generalization of   \cite[Lemma 2.2.5]{W}, relates the notion of  bimodule homomorphism and the one of a module generalized $d$-derivation when $d$ is inner. It answers the question concerning the converse implication of the one given in assertion (1) of Theorem   \ref{thm-inn-bimo}.

\begin{proposition}  \label{prop-gen-bimo}   Consider a module generalized $d$-derivation $ S: \mathcal{M}\longrightarrow \mathcal{M} $, where $ d  : \mathcal{A}\longrightarrow \mathcal{A} $  is a derivation. If
$d=[a_0,-]$ is inner, for some $ a_{0}\in \mathcal{A}$, then there is a  bimodule homomorphism $\Phi: \mathcal{M}\longrightarrow \mathcal{M} $ such that $S=\Phi + [a_0,-]$.\\
Moreover, $S$ is inner if and only if $\Phi$ is inner.
\end{proposition}
\pr  As done in \cite[Lemma 2.2.5]{W}, we simply need to prove that the map  $S-[a_0,-]$ is   a bimodule homomorphism.\cqfd

Note that for every $a_0\in Z(\mathcal{A})$ and every $b_0\in  \mathcal{A}$, the derivations $[a_0,-]$ and $[a_0+b_0,-]$ on $\mathcal{A}$ are the same. Thus, the fact that  $\Phi$ is a central  inner  bimodule homomorphism  will not affect the coresponding element of the inner derivation $d$ associated to the module generalized derivation $S$.\medskip

We end this section with a  characterization of inner derivations on trivial extension algebras.  It is  a generalization of \cite[Proposition 2.2.3]{W} which  characterizes inner derivations on   triangular matrix algebras (compare it also with  \cite[Proposition 2.2]{MP}).

\begin{theorem}\label{thm-form2}  
Consider a  derivation $ \mathcal{D}:\mathcal{A}\ltimes \mathcal{M}\longrightarrow \mathcal{A}\ltimes \mathcal{M}$ of the form
  $$    \mathcal{D} =(\mathcal{D}_{\mathcal{A}} +T ,\mathcal{D}_{\mathcal{M}} +S ), 
$$
  where $ \mathcal{D}_{\mathcal{A}},\;  T,\;  \mathcal{D}_{\mathcal{M}}$ and  $S$ as indicated in Lemma \ref{lem-der}. Then, the following assertions hold.
\begin{enumerate}
    \item If $\mathcal{D}$ is inner, then     $T=0$ and both $ S$ and $\mathcal{D}_{\mathcal{M}}$ are inner.
 \item The converse holds if we assume that $l.Ann_{\mathcal{A}}(\mathcal{M})\cap r.Ann_{\mathcal{A}}(\mathcal{M})=0$.
\end{enumerate}
\end{theorem}
\pr $\Rightarrow)$  Suppose that $\mathcal{D}$ is inner. Then  there exists $ (a_{0},m_{0})\in  \mathcal{A}\ltimes \mathcal{M}$ such that, for every $ a\in \mathcal{A}$ and $ m\in \mathcal{M}$,
\begin{eqnarray*}
  \mathcal{D}(a,m) &=& [(a_{0},m_{0}),(a,m)] \\
   &=& ([a_{0},a],[m_{0},a]+[a_{0},m]).
\end{eqnarray*}
 Then, for $a=0$, $D(0,m)=(T(m),S(m))=(0,[a_{0},m])$. This shows that $T(m)=0$ and $S(m)=[a_{0},m]$. Now, for $m=0$, we get $\mathcal{D}_{\mathcal{M}}(a)=[m_{0},a]$, as desired.\smallskip

 $\Leftarrow)$ Now we assume that $l.Ann_{\mathcal{A}}(\mathcal{M})\cap r.Ann_{\mathcal{A}}(\mathcal{M})=0$. Since $\mathcal{D}_{\mathcal{M}}$ and $S$ are inner, there are $a_{0}\in \mathcal{A}$ and $ m_{0}\in \mathcal{M}$ such that, for every $ a\in \mathcal{A}$ and $ m\in \mathcal{M}$, $\mathcal{D}_{\mathcal{M}}(a)=[ m_{0},a]$ and $S(m)=[a_{0},m]$. By  assertion (1) of Theorem  \ref{thm-inn-bimo},   $\mathcal{D}_{\mathcal{A}}(a)=[a_{0},a]$. Therefore,
 \begin{eqnarray*}
  \mathcal{D}(a,m) &=&([a_{0},a],[m_{0},a]+[a_{0},m])  \\
   &=& [(a_{0},m_{0}),(a,m)],
\end{eqnarray*}
as desired. \cqfd

It is worth noting that there are examples of trivial extension algebras showing that the converse implication in Theorem \ref{thm-form2} above does not hold true if we drop the condition $l.Ann_{\mathcal{A}}(\mathcal{M})\cap r.Ann_{\mathcal{A}}(\mathcal{M})=0$.\\

When $\mathcal{M}$ is assumed to be  $\mathcal{A}$, the condition of Theorem  \ref{thm-form2} is fulfilled. Thus, we get  the following corollary (compare it with \cite[Proposition 2.4]{MP}).

\begin{corollary} \label{cor-form2} Assume $\mathcal{M}=\mathcal{A}$ and consider a  derivation $ \mathcal{D}:\mathcal{A}\ltimes \mathcal{M}\longrightarrow \mathcal{A}\ltimes \mathcal{M}$ of the form
  $$    \mathcal{D} =(\mathcal{D}_{\mathcal{A}} +T ,\mathcal{D}_{\mathcal{M}} +S )
$$
  where $ \mathcal{D}_{\mathcal{A}},\;  T,\;  \mathcal{D}_{\mathcal{M}}$ and  $S$   as indicated in Lemma \ref{lem-der}. Then $\mathcal{D}$ is inner
     if and only if   $T=0$,  $\mathcal{D}_{\mathcal{M}}$ is inner, $ S= \mathcal{D}_{\mathcal{A}}  $ and $S$ is inner.\\
Precisely, the inner derivations on $\mathcal{A}\ltimes \mathcal{A}$ are only of the form $(d_1  ,d_0+d_1 )$, where $d_0$ and $d_1$ are   inner  derivations on $\mathcal{A}$.
\end{corollary}
\pr Only the assertion  $ S= \mathcal{D}_{\mathcal{A}}  $, when  $\mathcal{D}$ is inner, merits a proof. 
Thus, assume that  $\mathcal{D}$ is inner. Then,  by
Theorem \ref{thm-form2}, $S=[a_{0},-]$ is inner for some $a_0\in \mathcal{A}$. Then, the associated derivation $\mathcal{D}_{\mathcal{A}}=[a_{0},-]$ is also inner (by Theorem  \ref{thm-inn-bimo}). Now, by Proposition   \ref{prop-gen-bimo},    there is a  bimodule homomorphism $\Phi: \mathcal{M}\longrightarrow \mathcal{M} $ such that $S=\Phi + [a_0,-]$. Since   $S$ is inner,   $\Phi$ is also inner. Then, using assertion (2) of Theorem  \ref{thm-inn-bimo},    $a_0$ is in the  center $Z(\mathcal{A})$. This implies that $\Phi=0$, as desired.\cqfd


\section{First cohomology group of trivial extension algebras}

In this section  we study the first cohomology group of   trivial extension algebras. For this we start with the following observation. Let us note by $\der(\mathcal{A}\ltimes \mathcal{M})$  the subgroup of $\Der(\mathcal{A}\ltimes \mathcal{M})$ of derivations on $\mathcal{A}\ltimes \mathcal{M}$ of the form $(d,S)$, where $d$ is  a derivation on $\mathcal{A}$ and $S: \mathcal{M} \to \mathcal{M}$ is a module generalized $d$-derivation. Also  we use  $\innder(\mathcal{A}\ltimes \mathcal{M})$ to denote the inner derivations on $\mathcal{A}\ltimes \mathcal{M}$  of the form $[(a_0,0),-]=([a_0,-],[a_0,-])$. Note that   $\innder(\mathcal{A}\ltimes \mathcal{M})$ is a subgroup of $\der(\mathcal{A}\ltimes \mathcal{M})$. Thus, we may consider the quotient group  $\der(\mathcal{A}\ltimes \mathcal{M})/\innder(\mathcal{A}\ltimes \mathcal{M})$ which will be called, by analogy with the classical case, the \textit{restricted first cohomology group}  of $\mathcal{A}\ltimes \mathcal{M}$, and denoted by $h^{1}(\mathcal{A}\ltimes \mathcal{M})$.\\
Recall   the first cohomology group $ H^1(\mathcal{A},\mathcal{M}) = \Der(\mathcal{A},\mathcal{M}) /  \Innder(\mathcal{A},\mathcal{M})$. When $\mathcal{A}=\mathcal{M}$,  $ H^1(\mathcal{A},\mathcal{M})$ is simply denoted by $ H^1(\mathcal{A})$.

\begin{proposition}\label{Prop-coho}
The following assertions hold.
\begin{enumerate}
    \item There is a natural group homomorphism: $$\Der(\mathcal{A}\ltimes \mathcal{M})\cong  \Der(\mathcal{A},\mathcal{M})   \oplus \der(\mathcal{A}\ltimes \mathcal{M}) \oplus \mathcal{E}(\mathcal{M}) . $$
    \item There is a natural group homomorphism: $$\Innder(\mathcal{A}\ltimes \mathcal{M})\cong  \Innder(\mathcal{A},\mathcal{M})   \oplus \innder(\mathcal{A}\ltimes \mathcal{M}) . $$
\end{enumerate}
Consequently, there is  a natural group homomorphism:
 $$H^1(\mathcal{A}\ltimes \mathcal{M})\cong  H^1(\mathcal{A},\mathcal{M})   \oplus h^1(\mathcal{A}\ltimes \mathcal{M}) \oplus \mathcal{E}(\mathcal{M}). $$
\end{proposition}
 \pr The assertions are simple consequences of Lemma \ref{lem-der} and Theorem  \ref{thm-form2}.\cqfd

Now we give our main result which relates $h^1(\mathcal{A}\ltimes \mathcal{M})$, $H^1(
 \mathcal{A})$ and the quotient group  $\End_{\mathcal{A}-\mathcal{A} }(\mathcal{M})/\Innbi_{\mathcal{A} }(\mathcal{M}) $. By analogy with the classical notation, let us denote by $ H^1_{\mathcal{A}}(\mathcal{M})$ the quotient group      $\End_{\mathcal{A}-\mathcal{A} }(\mathcal{M})/\Innbi_{\mathcal{A} }(\mathcal{M})$. Using Proposition \ref{Prop-coho},  the main result provides a generalization of  the classical result \cite[Theorem 5.5]{CMRS} and the recent result  \cite[Theorem  4.4]{AGST} (see  Corollary \ref{cor-triv-princ}). To see this fact, we give a brief discussion of the scope. \medskip

In \cite{CMRS} the first cohomology group of a particular case of trivial extensions is studied.  Namely, by \cite[Theorem 5.5]{CMRS}, we have, if $\mathcal{R}$ is a field, $\mathcal{A}$ is assumed to be a finite dimensional algebra and $\mathcal{M}=D\mathcal{A}$ is the dual  $\mathcal{A}$-bimodule of $\mathcal{A}$,
 then  $$H^1(\mathcal{A}\ltimes \mathcal{M})\cong  H^1(\mathcal{A},\mathcal{M})   \oplus  H^1(\mathcal{A}) \oplus Z(\mathcal{A})   \oplus Alt_\mathcal{A}(D\mathcal{A}), $$ where
$Alt_\mathcal{A}(D\mathcal{A})$ is the set of skew-symmetric bilinear forms $\beta$ over $D\mathcal{A}$ such that
$\beta(fa, g)=\beta(f, ag)$   for all $f, g \in D\mathcal{A}$ and $a \in \mathcal{A}$.  As noted above \cite[Theorem 5.5]{CMRS}, this vector space coincides with $ \mathcal{E}(D\mathcal{A})$. Also, note that,  from \cite[Proposition 3.3 and Example 3.5]{CMRS},   the center $Z(\mathcal{A})$   of $\mathcal{A}$ has a symmetric action on $D\mathcal{A}$  (that is, $af=fa$ for every  $a\in Z(\mathcal{A})$  and  $f\in D\mathcal{A}$). This shows that   $ \Innbi_{\mathcal{A} }(\mathcal{M})=0$ and so $H^1_{\mathcal{A}}(\mathcal{M})=\End_{\mathcal{A}-\mathcal{A} }(\mathcal{M})$. On the other hand,
 $$\End_{\mathcal{A}-\mathcal{A} }(\mathcal{M})=\End_{\mathcal{A}-\mathcal{A} }(\mathcal{A})=Z(\mathcal{A}).$$
  Then, by \cite[Theorem 5.5]{CMRS}, we deduce that  $$h^{1}(\mathcal{A}\ltimes \mathcal{M})= H^1_{\mathcal{A}}(\mathcal{M}) \oplus  H^1(\mathcal{A}).$$
Thus, using \cite[Corollary  3.7]{AGST}, one can show that the main result provides a generalization of \cite[Theorem 5.5]{CMRS}. In fact, it is a generalization of the more general recent result \cite[Theorem  4.4]{AGST} (see Corollary \ref{cor-triv-princ}). \medskip

Let us denote by $[a]$ the equivalence class of an element $a$ in a given quotient group. We will use the canonical projection $\pi_{\mathcal{A}}$  of $\der(\mathcal{A}\ltimes \mathcal{M})$  onto $\Der(\mathcal{A})$ (i.e., $\pi_{\mathcal{A}}((d,S))=d$, where $d$ is a derivation on $\mathcal{A}$ and $S: \mathcal{M} \to \mathcal{M}$ is a module generalized $d$-derivation). Also, we will use the linear application $  \Phi:   \End_{\mathcal{A}-\mathcal{A} }(\mathcal{M}) \longrightarrow \der(\mathcal{A}\ltimes \mathcal{M})$ defined by  $  \Phi(S)=(0,S)$ for all $S\in \End_{\mathcal{A}-\mathcal{A} }(\mathcal{M})$.

\begin{theorem}\label{thm-princ}
  Consider the two maps
$ \widehat{\pi}_{\mathcal{A}}:  h^1(\mathcal{A}\ltimes \mathcal{M}) \longrightarrow H^{1}(\mathcal{A})$ given by $  \widehat{\pi}_{\mathcal{A}}([\mathcal{D}])=[\pi_{\mathcal{A}}(\mathcal{D})]$  and
 $\widehat{ \Phi}:    H^1_{\mathcal{A}}(\mathcal{M}) \longrightarrow h^1(\mathcal{A}\ltimes \mathcal{M}) $ given by   $ \widehat{ \Phi}([S]):=[\Phi(S)]$. The maps $ \widehat{\pi}_{\mathcal{A}}$ and  $\widehat{ \Phi}$ are well defined group homomorphisms.  \\
Moreover, we have the following exact sequence of group homomorphisms
  $$
    \xymatrix{
        0 \ar[r]  &   H^1_{\mathcal{A}}(\mathcal{M})   \ar[r]^-{\widehat{ \Phi} } & h^1(\mathcal{A}\ltimes \mathcal{M})\ar[r]^-{ \widehat{\pi}_{\mathcal{A}}       }  & H^{1}(\mathcal{A}).
         }$$
\end{theorem}
\pr   
We prove that $\widehat{\pi}_{\mathcal{A}}$ is well defined. First note that $ \pi_{\mathcal{A}}$ sends any inner derivation of $\mathcal{A}\ltimes \mathcal{M}$ to an inner derivation of  $\mathcal{A}$. Now consider $\mathcal{D}_{1}=(d_1,S_1),\mathcal{D}_{2}=(d_2,S_2)\in \der(\mathcal{A}\ltimes \mathcal{M})$ such that  $[\mathcal{D}_{1}]=[\mathcal{D}_{2}]$ in $h^1(\mathcal{A}\ltimes \mathcal{M})$. Then, there exists $(a_0,0 )\in \mathcal{A}\ltimes \mathcal{M}$ such that $ \mathcal{D}_{1}=\mathcal{D}_{2}+[(a_0, 0),-]$.   Then  $d_{1}=d_{2}+[a_0,-]$. Therefore,   $[\pi_{\mathcal{A}}(\mathcal{D}_{1})]=[\pi_{\mathcal{A}}(\mathcal{D}_{2})]$.\\
Now   we show that $\widehat{\Phi}$ is well defined. Also note that $\Phi$ sends a central inner bimodule homomorphism to an inner
derivation of  $\mathcal{A}\ltimes \mathcal{M}$. Consider $f,g\in\End_{\mathcal{A}-\mathcal{A} }(\mathcal{M})$ such that   $ [f]=[g]$ in
 $ H^1_{\mathcal{A}}(\mathcal{M})$. Then, there exists $a_0 \in Z(\mathcal{A})$ such that $ f=g+[a_0,-]$.
  Then,  $ (0,f)=(0,g)+[(a_0,0),-]$  in $\der(\mathcal{\mathcal{A}\ltimes \mathcal{M}})$. Therefore, $ [(0,f)]=[(0,g)]$ in $h^1(\mathcal{A}\ltimes \mathcal{M})$.\\
We prove that $\widehat{ \Phi}$ is injective. Consider $f \in\End_{\mathcal{A}-\mathcal{A} }(\mathcal{M})$ such that   $ \widehat{ \Phi}([f])=0$. Then, $[(0,f)] = 0$ in $ h^1(\mathcal{A}\ltimes \mathcal{M}) $. This shows that there is   $a_0\in \mathcal{A}$ such that $(0,f)=[(a_0,0),-]$. Thus, as a derivation on $\mathcal{A}$,  $[a_0 ,-]=0$. Therefore, $a_0\in Z(\mathcal{A})$, and so $f$ is a central inner bimodule homomorphism, as desired.\\
It remains to prove that   $ \Ker(\widehat{\pi}_{A})=\Im(\widehat{\Phi})$. First note that  $\widehat{\pi}_{A} \widehat{\Phi} =0$ and hence  $ \Im(\widehat{\Phi})\subseteq \Ker(\widehat{\pi}_{A})$. For the converse, consider  $ D =(d ,S ) \in \der(\mathcal{A}\ltimes \mathcal{M})$ such that $\widehat{\pi}_{A}([D])=0$. This means that   $ d$ is inner, hence $d=[a_0,-]$ for some $a_0\in \mathcal{A}$. Then, by Proposition \ref{prop-gen-bimo}, there is   $f \in \End_{\mathcal{A}-\mathcal{A} }(\mathcal{M})$ such that $S=f+ [a_0,-]$. Then, in $ h^1(\mathcal{A}\ltimes \mathcal{M})$, we get $$[D]=[(d,f+[a_0,-])]=[(0,f)+[(a_0,0),-]]=[(0,f)]=\widehat{ \Phi}([f]),  $$
as desired.\cqfd

From \cite[Proposition 4.8]{AGST}\footnote{The authors would like to thank Professor Rachel Taillefer (Université Blaise Pascal,  France) for discussing with us  \cite[Proposition 4.8]{AGST}.}, we have $ H^1(\mathcal{A}\ltimes\mathcal{M} ,\mathcal{M})=  H^1(\mathcal{A}  ,\mathcal{M})\oplus \End_{\mathcal{A}-\mathcal{A}}(\mathcal{M})$. Then,    there is a natural surjective group homomorphism $$ H^1(\mathcal{A}\ltimes\mathcal{M} ,\mathcal{M})  \longrightarrow H^1(\mathcal{A}  ,\mathcal{M})\oplus H^1_{\mathcal{A}}(\mathcal{M}).$$
Obviously, this homomorphism is in fact the identity when the center $Z(\mathcal{A})$   of $\mathcal{A}$ has a symmetric action on $\mathcal{M}$. Thus, the following result is a  generalization of   \cite[Theorem  4.4]{AGST} (compare it also  with \cite[Theorem 2.5]{MP}).

\begin{corollary}\label{cor-triv-princ}
There is an exact sequence of group homomorphisms
  $$
    \xymatrix{
        0 \ar[r]  &   H^1_{\mathcal{A}}(\mathcal{M})  \oplus  H^1(\mathcal{A},\mathcal{M})    \oplus \mathcal{E}(\mathcal{M})
   \ar[r] & H^1(\mathcal{A}\ltimes \mathcal{M})\ar[r]^-{ \widehat{\pi}_{\mathcal{A}}       }  & H^{1}(\mathcal{A}),
         }$$
where the first (injective)   homomorphism  is induced from $\widehat{ \Phi} $ and the natural isomomorphisms  of Proposition \ref{Prop-coho}.
\end{corollary}
 
As mentioned  in Remark \ref{rem-deriv},  both $H^1(\mathcal{A},\mathcal{M}) =0$ and $ \mathcal{E}(\mathcal{M})=0$  in the case of    triangular matrix algebras $ \mathcal{A}\ltimes \mathcal{M}=\tri(A,\mathcal{M},B)$ (where $\mathcal{A}=A\times B $). Thus, the following result is a generalization of \cite[Theorem 2.3.6]{W}.

\begin{corollary}\label{cor-matrix-princ} If $ \mathcal{A}\ltimes \mathcal{M}=\tri(A,\mathcal{M},B)$ (where $\mathcal{A}=A\times B $), then  there is an exact sequence of group homomorphisms
  $$
    \xymatrix{
        0 \ar[r]   &   H^1_{\mathcal{A}}(\mathcal{M})  \ar[r]^-{\ \widehat{ \Phi}}  
    & H^1(\tri(A,\mathcal{M},B)) \ar[r]^-{ \widehat{\pi}_{\mathcal{A}}}  & H^{1}(\mathcal{A}).
         }$$ 
\end{corollary}

  Also as a simple consequence of Theorem \ref{thm-princ}, we set the following generalization of  \cite[Corollary 2.3.8]{W}.

\begin{corollary}\label{cor-princ-identifica} 
If $H^{1}(\mathcal{A})=0$, then we get the following group isomorphisms: 
  $$
\begin{array}{ccc}
         H^1_{\mathcal{A}}(\mathcal{M})    \cong h^1(\mathcal{A}\ltimes \mathcal{M}) & \textrm{and} &   H^1_{\mathcal{A}}(\mathcal{M})  \oplus  H^1(\mathcal{A},\mathcal{M})    \oplus \mathcal{E}(\mathcal{M})
   \cong H^1(\mathcal{A}\ltimes \mathcal{M}).
       \end{array}  $$
\end{corollary}

 Naturaly one can ask of a relation between    $\End_{\mathcal{A}-\mathcal{A} }(\mathcal{M})/\InnBi_{\mathcal{A} }(\mathcal{M})$,  $h^1(\mathcal{A}\ltimes \mathcal{M})$ and $ H^{1}(\mathcal{A})$ as done in Thereom \ref{thm-princ}. Using, Theorem \ref{thm-inn-bimo},   we deduce that  $\End_{\mathcal{A}-\mathcal{A} }(\mathcal{M})/\InnBi_{\mathcal{A} }(\mathcal{M})=H^1_{\mathcal{A}}(\mathcal{M})$ when $l.Ann_{\mathcal{A}}(\mathcal{M})\cap r.Ann_{\mathcal{A}}(\mathcal{M})=0$. Then we get the following result.

\begin{corollary}\label{cor-princ-h1=H1}
 If $l.Ann_{\mathcal{A}}(\mathcal{M})\cap r.Ann_{\mathcal{A}}(\mathcal{M})=0$, 
  then we have the following exact sequence of group homomorphisms
  $$
    \xymatrix{
        0 \ar[r]  &   h^1_{\mathcal{A}}(\mathcal{M})   \ar[r]^-{\widehat{ \Phi} } & h^1(\mathcal{A}\ltimes \mathcal{M})\ar[r]^-{ \widehat{\pi}_{\mathcal{A}}       }  & H^{1}(\mathcal{A}),
         }$$
where   $h^1_{\mathcal{A}}(\mathcal{M})$ denotes the quotient group     $\End_{\mathcal{A}-\mathcal{A}}(\mathcal{M})/\InnBi_{\mathcal{A} }(\mathcal{M})$.
\end{corollary}

The following example shows that the result of Corollary \ref{cor-princ-h1=H1}  does not hold true without the condition $l.Ann_{\mathcal{A}}(\mathcal{M})\cap r.Ann_{\mathcal{A}}(\mathcal{M})=0$.

\begin{example}\label{exm-All-inn-1}
Consider $\mathcal{R}$, $A$, $B$, $M$,  $\mathcal{S}$ and $\mathcal{N}$ as in example  \ref{exm-inn-bimo}. Then, using Example \ref{exm-inn-bimo} and Remark \ref{rem-deriv}, we have $ H^{1}(\mathcal{S})=0$, $h^1_{\mathcal{S}}(\mathcal{N})=0$,
 $\mathcal{E}( \mathcal{N})=0$ and $H^1(\mathcal{S},\mathcal{N})=0$. However, $h^1(\mathcal{S}\ltimes \mathcal{N})\neq 0$ (so does $H^1(\mathcal{S}\ltimes \mathcal{N})$). Indeed,   consider the two elements $a_0 =                                                                                                                                    \begin{pmatrix}                                                                                                                                      1 & 0\\                                                                                                                                0 & 0                                                                                                                              \end{pmatrix}  $ and $b_0 =                                                                                                                                       \begin{pmatrix}                                                                                                                                       0 & 0\\                                                                                                                                0 & 1                                                                                                                               \end{pmatrix}                                                                                                                               $
of $\mathcal{S}$. Consider the map $D: \mathcal{S}\ltimes \mathcal{N}\to \mathcal{S}\ltimes \mathcal{N} $ defined as follows:
For every  $x=\begin{pmatrix}                                                                                                                                       a & m\\                                                                                                                                0 & b                                                                                                                               \end{pmatrix}  $ in $\mathcal{S}$ and $n\in  \mathcal{N} $,
$$D(x,n)=([a_0,x], [b_0,n])=(\begin{pmatrix}                                                                                                                                       0 & m\\                                                                                                                                0 & 0                                                                                                                               \end{pmatrix} ,n). $$
From Example \ref{exm-inn-bimo}, $[b_0,-]$ is an $[a_0,-]$-generalized derivation. Then, $D$ is a derivation on $\mathcal{S}\ltimes \mathcal{N}$. However, one can show that  $D$  is not inner. Indeed,  suppose that $D$ is inner. Then, there are   $c=\begin{pmatrix}                                                                                                                                       \alpha & e\\                                                                                                                                0 & \beta                                                                                                                               \end{pmatrix}  $ in $\mathcal{S}$ and $n_0\in \mathcal{N}$ such that $D=[(c,n_0),-]$. Then, after a simple calculation, using the module actions defined in Example \ref{exm-principal}, we deduce that, for every $m,n\in  M=\mathcal{N}=\mathcal{R}$, $\alpha m - m \beta =m$ and $\beta n - n\alpha=n$. This is impossible, as desired.
\end{example}

As an important consequence of Corollary \ref{cor-triv-princ}, we get a result that studies the innerness of derivations on $ \mathcal{A}\ltimes \mathcal{M}$.  \\

We use $ \overline{\Gamma}$ to denote the subgroup $ \pi_{\mathcal{A}}( \der(\mathcal{A}\ltimes \mathcal{M}) )$ of $\Der(\mathcal{A})$.

\begin{corollary}\label{cor-princ-inner-triv}
All  derivations  on  $\mathcal{A}\ltimes \mathcal{M}$ are inner if and only if the following assertions hold.
 \begin{enumerate}
    \item Every derivation  in $  \overline{\Gamma} $ is inner.
     \item  Every derivation  in   $\Der(\mathcal{ A},\mathcal{M})$ is inner.
    \item  Every $\mathcal{A}$-bimodule homomorphism  $\End_{\mathcal{A}-\mathcal{A} }(\mathcal{M})$ is  central inner.
  \item   $\mathcal{E}( \mathcal{M})=0$.
\end{enumerate}
\end{corollary}

Note that when $ \mathcal{A}\ltimes \mathcal{M}=\tri(A,\mathcal{M},B)$ (where $\mathcal{A}=A\times B)$ a derivation $d$ of $A\times B$ is of the form $d=(d_1,d_2)$ where $d_1$ and $d_2$ are derivations of $A$ and $B$ respectively. Then,   a module generalized $(d_1,d_2)$-derivation $ S: \mathcal{M}\longrightarrow \mathcal{M} $ is of the form  $ S(am)=aS(m)+d(a)m=aS(m)+d_1(a)m$ and $S(ma)=S(m)a+md(a)=S(m)a+md_2(a)$
for  all $m\in \mathcal{M} $ and $a\in \mathcal{A}$. Thus,  $ \overline{\Gamma}$ used in Corollary \ref{cor-princ-inner-triv}  is naturally a group  product  $ \overline{\Gamma}= \overline{\Gamma}_1 \times  \overline{\Gamma}_2$  where $\overline{\Gamma}_1$ and $\overline{\Gamma}_2$  are constitued with the associated derivations on $A$ and $B$ respectively. With these notations and using Corollary \ref{cor-princ-inner-triv} and  Remark \ref{rem-deriv}, we   recover  Benkovi\v{c}'s result \cite[Theorem 3.5]{D}.

\begin{corollary} \label{cor-princ-inner-matr}
 Let $A$ and $B$ be unital algebras over a commutative ring $C$  and let $\mathcal{M}$ be a unital $(A,B)$-bimodule. Then all the derivations  on   $\tri(A,\mathcal{M},B)$ are inner if and only if the
 following assertions hold.
 \begin{enumerate}
    \item  Every derivation   in $  \overline{\Gamma}_1 $ is inner.
    \item   Every derivation   in $ \overline{\Gamma}_2 $ is inner.
\item  Every  bimodule homomorphism  in $\End_{A\times B}(\mathcal{M})$ is central inner.
\end{enumerate}
\end{corollary}

Triangular matrix algebras can be used as examples of trivial extension algebras  satisfying the conditions of  Corollary \ref{cor-princ-inner-triv}. Take, for instance, the algebra $\tri(\R,\R,\R)$. However, naturally, one may ask for an example of a trivial extension algebra which   has not a   triangular matrix representation and   satisfies the conditions of  Corollary \ref{cor-princ-inner-triv}. The authors have not been able to give such an example. In fact, every   studied example   has a triangular matrix representation.  This leads  to the following natural question.\medskip

\noindent\textbf{Question.} Does the condition  $H^1(\mathcal{A}\ltimes \mathcal{M})=0$ imply that $\mathcal{A}\ltimes \mathcal{M}$ has a   triangular matrix representation?\medskip

Observations of some studied examples show that the key for reaching this target could be the study of the property  $(3)$ in  Corollary \ref{cor-princ-inner-triv}.  We end this paper by showing that this property  with some other mild conditions  assure that the trivial extension  has a  triangular matrix representation. This gives  a partial affirmative answer to the question above. \medskip

The following lemma gives an evident, but important, consequence of the property  $(3)$ in  Corollary \ref{cor-princ-inner-triv}.

\begin{lemma}\label{lem-prin-1} If every bimodule homomorphism on $\mathcal{M}$ is central inner, then for every $a\in Z(\mathcal{A})$, there exists   $l_a\in Z(\mathcal{A})$ (resp.   $r_a\in Z(\mathcal{A}$)) such that, for all $m\in \mathcal{M}$,  $am=l_a m-ml_a$ (resp. $ma=r_a m-mr_a$).\\
In particular, there exists   $c\in Z(\mathcal{A})$ such that  $m=cm-mc$ for all $m\in \mathcal{M}$.
\end{lemma}
\pr  The result follows since   $\varphi_r: \mathcal{M} \longrightarrow \mathcal{M}$ defined by $\varphi_{r}(m)=ma$ and $\varphi_l: \mathcal{M} \longrightarrow \mathcal{M}$ defined by $\varphi_{l}(m)=am$  are bimodule homomorphisms.\cqfd

Note that, if the element  $c$  in  Lemma \ref{lem-prin-1}   is a zero-divisor, then $x\mathcal{M}x=0$ for every $x\in  \mathcal{A}$ that annihilates $c$ (either  on the left or on the right).\medskip

\begin{theorem}\label{thm-triv-is-trian}
  Assume   that  $l.Ann_\mathcal{A}(\mathcal{M})\cap r.Ann_\mathcal{A}(\mathcal{M})=0$.\\
 The trivial extension algebra $\mathcal{A}\ltimes \mathcal{M}$ has a   triangular matrix representation if and only if  there exists   $ c\in r.Ann(\mathcal{M})$ (or $c\in l.Ann(\mathcal{M})$) such that  $m=cm-mc$ for all $m\in \mathcal{M}$.
\end{theorem}
\pr  One implication is a consequence of the discussion in Section 2.\\  
We prove the ``if" part. We may suppose that $ c\in r.Ann(\mathcal{M})$. Similarly we prove the result if  $c\in l.Ann(\mathcal{M})$. \\
Consider $m\in  \mathcal{M}$. Since  $ c\in r.Ann(\mathcal{M})$, $m(c-c^2)=0$, and by Lemma  \ref{lem-prin-1}, $(c-c^2)m=c(1-c)m=-cmc=0$. Then, the fact that  $l.Ann_\mathcal{A}(\mathcal{M})\cap r.Ann_\mathcal{A}(\mathcal{M})=0$ implies that $c=c^2$.  Using Proposition \ref{prop-triv-trian}, it remains to show that  $(1-c)\mathcal{A}c=0$. Let    $m\in \mathcal{M}$ and let  $a\in \mathcal{A}$. Since  $ acm \in \mathcal{M}$,  $(1 - c)acm=-acmc = 0$. Also, since $m(1 - c)a \in \mathcal{M}$, $ m(1 - c)ac = 0$. Therefore,  $(1 - c)ac\in  l.Ann_\mathcal{A}(\mathcal{M})\cap r.Ann_\mathcal{A}(\mathcal{M})=0$.\cqfd

\noindent{\bf Acknowledgements.} The authors are deeply grateful to the referee for careful reading of the manuscript and helpful suggestions.


\end{document}